\documentclass[11pt,a4paper,reqno]{amsart}
\usepackage{mathrsfs,amssymb}

\newcommand{\scrA}{\mathscr{A}}
\newcommand{\scrB}{\mathscr{B}}
\newcommand{\scrR}{\mathscr{R}}
\newcommand{\scrS}{\mathscr{S}}

\newcommand{\presar}{\langle \scrA \mid \scrR \rangle}
\newcommand{\presas}{\langle \scrA \mid \scrS \rangle}
\newcommand{\presbs}{\langle \scrB \mid \scrS \rangle}

\newcommand{\ol}[1]{\overline{#1}}
\newcommand{\nequiv}{\not\equiv}

\newtheorem{theorem}{Theorem}[section]

\newtheorem{corollary}[theorem]{Corollary}

\newtheorem{proposition}[theorem]{Proposition}
\newtheorem{lemma}[theorem]{Lemma}

\newcounter{fiddletheoremtemp}

\begin{document}

\title[Decision Problems and Properties of Small Overlap Monoids]{Uniform Decision Problems and Abstract Properties of Small Overlap Monoids}


\maketitle

\begin{center}

    MARK KAMBITES

    \medskip

    School of Mathematics, \ University of Manchester, \\
    Manchester M13 9PL, \ England.

    \medskip

    \texttt{Mark.Kambites@manchester.ac.uk} \\

    \medskip

\end{center}

\begin{abstract}
We study the way in which the abstract structure of a small overlap monoid
is reflected in, and may be algorithmically deduced from, a small
overlap presentation. We show that every $C(2)$ monoid admits an essentially
canonical $C(2)$ presentation; by counting canonical presentations we
obtain asymptotic estimates for the number of non-isomorphic monoids
admitting $a$-generator, $k$-relation presentations of a given length.
We demonstrate
an algorithm to transform an arbitrary presentation for a $C(m)$ monoid
($m$ at least $2$) into a canonical
$C(m)$ presentation, and a solution to the isomorphism
problem for $C(2)$ presentations.
We also find a simple combinatorial condition on a $C(4)$ presentation
which is necessary and sufficient for the monoid presented to be
left cancellative. We apply this to obtain algorithms to decide if a given
$C(4)$ monoid is left cancellative, right cancellative or cancellative,
and to show that cancellativity properties are \textit{asymptotically
visible} in the sense of generic-case complexity.
\end{abstract}

\section{Introduction}

Small overlap conditions are natural combinatorial conditions on semigroup and monoid
presentations, which serve to limit the complexity of derivation sequences
between equivalent words. First studied by Remmers \cite{Remmers71,Remmers80}
in the 1970's,
they are the natural semigroup-theoretic analogue of the \textit{small cancellation conditions}
extensively used in combinatorial and geometric group theory \cite{Lyndon77}.
Like word hyperbolicity for groups \cite{Olshanskii92}, small overlap properties
are \textit{generic} for monoid presentations with a fixed number of generators and relations, in
the sense that a randomly chosen $a$-generator, $k$-relation presentation
of size $n$ satisfies any given small overlap condition with probability
which approaches 1 as $n$ increases \cite{K_generic}.

Remmers' original work showed that monoids with presentations satisfying the
condition $C(3)$ have decidable word problem; an accessible account of this
and related results can be found in \cite{Higgins92}. Recent research of the
author \cite{K_solupdate,K_smallover1,K_smallover2} has shown that monoids with presentations
satisfying the slightly stronger condition $C(4)$ have linear time solvable
word problem, and a regular language of linear-time computable normal forms,
and are also rational (in the sense of Sakarovitch \cite{Sakarovitch87}), asynchronous
automatic and word hyperbolic (in the sense of Duncan and Gilman
\cite{Duncan04}).

This paper is devoted to the question of how abstract
algebraic properties (by which we mean isomorphism invariants) of small
overlap monoids are reflected in, and can be
algorithmically deduced from, their presentations.
Section~\ref{sec_preliminaries} briefly recalls some key definitions
and results from the theory of small overlap presentations.
In Section~\ref{sec_minimal} we consider isomorphisms between small
overlap monoids. We show that every isomorphism type of small overlap
monoids has a canonical
small overlap presentation, and that this presentation
can be effectively computed from any other presentation for the monoid.
It follows that the \textit{isomorphism problem} for small overlap monoids is
decidable. By counting canonical presentations, we obtain asymptotic
estimates for the number of isomorphism types of $a$-generator $k$-relation
semigroups of a given size, as a function of the positive integers $a$ and
$k$.

In Section~\ref{sec_canc} we turn our attention to cancellativity properties.
By applying results from \cite{K_solupdate,K_smallover1},
we are able to give elementary combinatorial characterisations of those $C(4)$
presentations which present left cancellative, right cancellative and
cancellative monoids. Since these properties of presentations can be easily
tested, it follows that one can efficiently decide whether a given $C(4)$
presentation presents a left, right or two-sided cancellative monoid.
We deduce also that left cancellativity, right cancellativity and cancellativity
are \textit{asymptotically
visible} properties of finite presentations, in the sense that the
proportion of $\scrA$-generated, $k$-relation monoid presentations of size $n$
which present monoids with any of these properties converges to a limit
strictly between $0$ and $1$, as $n$ increases.

\section{Small Overlap Monoids}\label{sec_preliminaries}

In this section we recap the definitions and some key results
concerning small overlap monoids.
We begin by recalling the basic definitions of combinatorial semigroup
theory, chiefly in order to fix notation and terminology.
Let $\scrA$ be an \textit{alphabet}, or set of symbols. A \textit{word} over
$\scrA$ is a finite sequence of zero or more symbols from $\scrA$. The
\textit{free monoid} $\scrA^*$ is the set of all words over $\scrA$, with multiplication
defined by concatenation of sequences. The unique \textit{empty word}
containing no symbols is denoted $\epsilon$; it forms the identity element
in $\scrA^*$. A \textit{monoid presentation} consists of a pair $\presar$
where $\scrA$ is an alphabet, and $\scrR$ is a binary relation on the
free monoid $\scrA^*$. The presentation is called \textit{finite} if both
$\scrA$ and $\scrR$ are finite. The elements of the binary relation $\scrR$,
which are pairs from $\scrA^* \times \scrA^*$, are rather confusingly called
the \textit{relations} of the presentation. A relation with both sides the
same (that is, one of the form $(u,u)$ for some $u \in \scrA^*$) is termed
\textit{trivial}.

The \textit{maximum relation length} of the presentation is the length of
the longest word appearing as one side of a relation in $\scrR$ (or $0$ if
there are no relations), while the \textit{sum relation length} is the
total length of all the words forming sides of relations in $\scrR$.

We say that a
word $u \in \scrA^*$ is obtained from a word $v \in \scrA^*$ by
\textit{an application of a relation from $\scrR$} if $u = pxq$ and
$v = pyq$ for some words $p,q,x,y \in \scrA^*$ such that $(x,y) \in \scrR$
or $(y,x) \in \scrR$. We define a binary relation $\equiv_\scrR$ (called
just $\equiv$ where there is no ambiguity as to the presentation) on $\scrA^*$
by $u \equiv_\scrR v$ if $v$ can be obtained by $u$ by a finite sequence of
zero or more applications of relations from $\scrR$. In fact $\equiv_\scrR$
is a \textit{congruence} on $\scrA^*$, that is, an equivalence relation
compatible with the multiplication. We denote by $[u]_\scrR$
(or just $[u]$) the equivalence class of the word $u \in \scrA^*$. The
equivalence classes form a monoid with multiplication well-defined by
$$[u]_\scrR [v]_\scrR = [uv]_\scrR;$$
this is called the \textit{monoid presented by the presentation}.

We say that a word $p$ is a \textit{possible prefix} of $u$ if there exists a
(possibly empty) word $w$ with $pw \equiv u$, that is, if the element
represented by $u$ lies in the right ideal generated by the element
represented by $p$. The empty word is denoted $\epsilon$.

A \textit{relation word} is a word which occurs as one side of a
relation in the presentation. A \textit{piece} is a word in the
generators which occurs as a factor in sides of two \textit{distinct} relation
words, or in two different (possibly overlapping) places within one
side of a relation word.
 By convention, the empty word is always a piece.
We say that a presentation is \textit{$C(n)$}, where $n$ is a positive
integer, if no relation word can be written as the product of \textit{strictly
fewer than} $n$ pieces. Thus for each $n$, $C(n+1)$ is
a stronger condition than $C(n)$.

Notice that it is permissible for the same relation word to appear twice
or more in a $C(n)$ presentation, since by doing so it does not become a
factor of two \textit{distinct} relation words. We say that a presentation
is \textit{strongly} $C(n)$ if it is $C(n)$ and additionally has no repeated
relation words.
The condition we have called $C(n)$ is that used in \cite{Higgins92,Remmers71,Remmers80}
and was called \textit{weakly} $C(n)$ in \cite{K_solupdate}, while the condition
we have called \textit{strongly} $C(n)$ was called $C(n)$ in
\cite{K_smallover1,K_smallover2,K_generic}.
See \cite{K_solupdate} for detailed explanation of the relationship
between the two definitions.

We say that a relation word $\ol{R}$ is a \textit{complement} of a relation
word $R$ if there are relation words $R = R_1, R_2, \dots, R_n = \ol{R}$ such
that either $(R_i, R_{i+1})$ or $(R_{i+1}, R_i)$ is a relation in the
presentation for $1 \leq i < n$. We say that $\ol{R}$ is a \textit{proper}
complement of $R$ if, in addition, $\ol{R} \neq R$. Abusing notation and
terminology slightly, if $R = X_R Y_R Z_R$ and
$\ol{R} = X_{\ol{R}} Y_{\ol{R}} Z_{\ol{R}}$ then we
write $\ol{X_R} = X_{\ol{R}}$, $\ol{X_R Y_R} = X_{\ol{R}} Y_{\ol{R}}$ and so
forth. We say that $\ol{X_R}$ is a complement of $X_R$, and $\ol{X_R Y_R}$
is a complement of $X_R Y_R$.

Now let $\presar$ be a $C(3)$ presentation. Recall that a \textit{relation prefix} of a word
is a prefix which admits a (necessarily unique, as a consequence of the
small overlap condition) factorisation of the form $a X Y$ where $X$ and $Y$
are the maximal piece prefix and middle word respectively of some relation
word $XYZ$. An \textit{overlap prefix (of length $n$)} of
a word $u$ is a relation prefix which admits an (again necessarily unique)
factorisation of the form $b X_1 Y_1' X_2 Y_2' \dots X_n Y_n$ where
\begin{itemize}
\item $n \geq 1$;
\item $b X_1 Y_1' X_2 Y_2' \dots X_n Y_n$ has no factor of the form $X_0Y_0$,
where $X_0$ and $Y_0$ are the maximal piece prefix and middle word respectively
of some relation word, beginning before the end of the prefix $b$;
\item for each $1 \leq i \leq n$, $R_i = X_i Y_i Z_i$ is a relation word with
$X_i$ and $Z_i$ the maximal piece prefix and suffix respectively; and
\item for each $1 \leq i < n$, $Y_i'$ is a proper, non-empty prefix of $Y_i$.
\end{itemize}
Notice that if a word has a relation prefix, then the shortest such must
be an overlap prefix. A relation prefix $a XY$ of a word $u$ is called
 \textit{clean} if $u$ does \textit{not} have a prefix
$$a XY' X_1 Y_1$$
where $X_1$ and $Y_1$ are the maximal piece prefix and middle word respectively
of some relation word, and $Y'$ is a proper, non-empty prefix of $Y$.

We recall some key technical results about $C(n)$ presentations.

\begin{proposition}[\cite{K_solupdate} Proposition~3]\label{prop_dumpprefix}
Suppose a word $u$ has an overlap prefix $w X Y$ and that
$u = w X Y u''$. Then $u \equiv v$ if and only if $v = w v'$ where
$v' \equiv X Y u''$.
\end{proposition}

\begin{lemma}[\cite{K_solupdate} Lemma~3]\label{lemma_eq}
Let $\presar$ be a $C(4)$ presentation.
Suppose $u = X Y u'$ where $XY$ is a clean overlap prefix of
$u$. Then $u \equiv v$ if and only if one of the following mutually
exclusive conditions holds:
\begin{itemize}
\item[(1)] $u = XYZ u''$ and $v = XYZ v''$ and $\ol{Z} u'' \equiv \ol{Z} v''$
for some complement $\ol{Z}$ of $Z$;
\item[(2)] $u = X Y u'$, $v = X Y v'$, and $Z$ fails to be a
prefix of at least one of $u'$ and $v'$, and $u' \equiv v'$;
\item[(3)] $u = X Y Z u''$, $v = \ol{X} \ol{Y} \ol{Z} v''$ for some
uniquely determined proper complement $\ol{XYZ}$ of $XYZ$,
and $\hat{Z} u'' \equiv \hat{Z} v''$ for some complement $\hat{Z}$
of $Z$;
\item[(4)] $u = X Y u'$, $v = \ol{X} \ol{Y} \ol{Z} v''$ for some uniquely
determined proper complement $\ol{XYZ}$ of $XYZ$ but
$Z$ is not a prefix of $u'$ and $u' \equiv Z v''$;
\item[(5)] $u = X Y Z u''$, $v = \ol{X} \ol{Y} v'$ for some uniquely
determined
proper complement
$\ol{XYZ}$ of $XYZ$,
but $\ol{Z}$ is not a prefix of $v'$ and $\ol{Z} u'' \equiv v'$;
\item[(6)] $u = X Y u'$, $v = \ol{X} \ol{Y} v'$ for some uniquely
determined proper complement
$\ol{XYZ}$ of $XYZ$, $Z$ is not
a prefix of $u'$ and $\ol{Z}$ is not a prefix of $v'$, but
$Z = z_1 z$, $\ol{Z} = z_2 z$, $u' = z_1 u''$, $v' = z_2 v''$ where
$u'' \equiv v''$ and $z$ is the maximal common suffix of $Z$ and $\ol{Z}$,
$z$ is non-empty, and $z$ is a possible prefix of $u''$.
\end{itemize}
\end{lemma}

\begin{proposition}[\cite{K_solupdate} Corollary 1]\label{prop_nomopnorel}
Let $\presar$ be a $C(3)$ presentation. If a word $u$ has no clean overlap
prefix, then it contains no relation word as a factor, and so if $u \equiv v$
then $u = v$.
\end{proposition}

\section{Minimal Presentations and Isomorphisms}\label{sec_minimal}

Our main aim in this section is to show that every small overlap
monoid has a \textit{canonical} small overlap presentation, and that
this presentation can be effectively computed from any other presentation
for the monoid. From this we are able to demonstrate a solution to the
isomorphism problem for monoid presentations satisfying the condition
$C(2)$. Some of the essential ideas in this section were prefigured in
work of Jackson \cite{Jackson91}, although the main results and
applications given here are new.

We begin by introducing some terminology.
A presentation $\presar$ is called an \textit{equivalence presentation} if the set
$\scrR$ of relations is an equivalence relation on the set of relation words.
The \textit{equivalence
closure} of a presentation $\presar$ is the presentation $\presas$ where
$\scrS$ is the reflexive, symmetric, transitive closure of $\mathscr{R}$.
Clearly, $\presas$
presents the same monoid as $\presar$.	
Notice also that, since we require the set of relations to be an equivalence
relation only on the set of relation words, and not on the whole free monoid $\scrA^*$,
the equivalence closure of a finite presentation is still finite, and can be
easily computed.
It has the same set of relation words as the original presentation. In
particular, it satisfies $C(m)$ for any $m \geq 1$ if and only if the original
presentation satisfies $C(m)$.

Recall that if $\langle \mathscr{A} \mid \mathscr{R} \rangle$ is a monoid
presentation then a generator $a \in \scrA$ is called \textit{redundant}
if it is equivalent to a product of (zero or more) other generators, or
equivalently, if the monoid presented by the presentation is generated by
the subset $\mathscr{A} \setminus \lbrace a \rbrace$.
We call a presentation \textit{generator-minimal} if it has no redundant
generators. Notice in particular that, in a generator-minimal presentation,
no two generators represent the same element. A non-identity
element $s$ of the monoid presented is called \textit{indecomposable} (or by some
authors, an \textit{atom}) if whenever $x,y$ are elements of the monoid such
that $xy = s$ we have $x=1$ or $y = 1$. A generator $a \in \scrA$ is called
indecomposable if $[a]_\scrR$ is indecomposable. Notice that a non-identity
indecomposable element of a monoid must belong to every generating set
for the monoid.

\begin{proposition}\label{prop_relationwordclass}
Let $\presar$ be a $C(2)$ presentation and $u$ a relation word. Then
for any $v \in \scrA^*$ we have $u \equiv v$ if and only if there exists
$n \geq 1$ and words $u = u_1, u_2, \dots, u_n = v$ such that for $1 \leq i < n$ either
$(u_i, u_{i+1}) \in \scrR$ or $(u_{i+1}, u_i) \in \scrR$.
\end{proposition}
\begin{proof}
One implication is immediate from the definitions. For the converse,
observe that since the presentation satisfies $C(2)$, no relation word
contains another relation word as a proper factor. It follows that the
only relations which can be applied to the word $u$ are relations of
the form $(u, u')$ or $(u', u)$. But now $u'$ is also a relation word, so
a simple inductive argument shows that any rewriting sequence taking $u$
to $v$ must have the given form.
\end{proof}

\begin{proposition}\label{prop_relationwordfactorclass}
Let $\presar$ be a $C(2)$ presentation and $u$ a relation word. Then
for any proper factor $u'$ of $u$ and any $v \in \scrA^*$ we have
$u' \equiv v$ if and only if $u'=v$.
\end{proposition}
\begin{proof}
Since the presentation satisfies $C(2)$, no relation word contains another
relation word as a proper factor. It follows that $u'$ contains no relation
word as a factor, and so no relation can be applied to it.
\end{proof}

\begin{proposition}\label{prop_indecomp}
Let $\presar$ be a monoid presentation satisfying $C(2)$, and let
$a \in \scrA$. Then $a$ is either indecomposable or redundant.
\end{proposition}
\begin{proof}
Since the presentation satisfies $C(1)$, the empty word
is not a relation word. If the presentation contains no non-trivial relations of the
form $(w, a)$ or $(a,w)$ then no relations are applicable to $a$, and so $a$ is
indecomposable. On the other hand, suppose the presentation does contain
a non-trivial relation of the form $(w,a)$ or $(a,w)$. Then $a$ is a relation word,
which since the presentation satisfies $C(2)$ means that $a$ is not a piece,
and so $a$ cannot feature in the word $w$. Thus, the given relation can
be used to rewrite $a$ as a product of the other generators, which means
that $a$ is redundant.
\end{proof}

\begin{proposition}\label{prop_genminexists}
Let $m \geq 2$. Every $C(m)$ monoid has a generator-minimal
$C(m)$ equivalence presentation, which can be effectively computed starting from any
$C(m)$ presentation for the monoid.
\end{proposition}
\begin{proof}
Let $\presar$ be a $C(m)$ presentation for $M$, and suppose that this
presentation is not generator-minimal. Let $a \in \mathscr{A}$ be a redundant generator.
Clearly some non-trivial relation
in $\scrR$ can be applied to the word $a$; since the presentation satisfies $C(1)$ neither
side of this relation can be the empty word, so this relation must have the
form $(a,w)$ or $(w,a)$ for some $w \in \mathscr{A}^+$ with $w \neq a$.

Let $\hat{a}$ denote the set of all words $w$ such that $(a, w)$ or $(w,a)$
is a relation and $w \neq a$. Since the presentation satisfies $C(2)$ and
$a$ is a relation word, $a$ cannot be a factor of any $w \in \hat{a}$. Let
$$\scrB = \scrA \setminus \lbrace a \rbrace$$
and
$$\scrS = \left( \scrR \setminus \lbrace (a, w), (w, a) \mid w \in \scrA^* \rbrace \right) \cup \lbrace (u,v) \mid u, v \in \hat{a} \rbrace.$$

Consider now the presentation $\presbs$. Since every relation word in
$\presbs$ is also a relation word in $\presar$, it is clear that $\presbs$
is a $C(m)$ presentation with strictly fewer generators than $\presar$.
We claim that $\presbs$ is a presentation for the same monoid as $\presar$,
Indeed, let $\iota : \scrB \to \scrA$ be the inclusion map. Then $\iota$
extends to a morphism $\hat\sigma : \scrB^* \to \scrA^*$ of free monoids, and since every
relation in $\scrS$ is satisfied in $\presar$, this induces a morphism
$\sigma : \presbs \to \presar$ well-defined by $\sigma([w]_\scrS) = [\hat\sigma(w)]_\scrR$. Now the image of $\sigma$ contains
$[w]_\scrR = [a]_\scrR$ and also $[b]_\scrR$ for every $b \in \scrB = \scrA \setminus \lbrace a \rbrace$,
so $\sigma$ must be surjective. Now suppose $u, v \in \scrB^*$ are such that
$\sigma([u]_\scrS) = \sigma([v]_\scrS)$. Then $\hat\sigma(u) \equiv_\scrR \hat\sigma(v)$, so there is a sequence of words
$$u = \hat\sigma(u) = y_0,\  y_1, \ \dots, \ y_n = \hat\sigma(v) = v \ \in \ \scrA^*$$
such that each $y_{i+1}$ can be obtained from $y_i$ by an application of
a single relation in $\scrR$. For each $i$, let $z_i$ be the word obtained
from $y_i$ by replacing any occurences of the letter $a$ with the word
$w$. Since $a$ is not a proper factor of any relation word in $\presar$,
it follows from the definition of $\scrS$ that each $y_{i+1}$ can be obtained
from $y_i$ by an application of a single relation from $\scrS$. But
$y_0 = u$ and $y_n = v$, so we deduce that $[u]_\scrS = [v]_\scrS$. Thus,
$\sigma$ is injective, and so is an isomorphism.

We have shown that given a non-generator-minimal $C(m)$ presentation, there always
exists a $C(m)$ presentation for the same monoid with strictly fewer
generators; it follows that a $C(m)$ presentation with the fewest possible
number of generators must be generator-minimal. Moreover, given a non-generator-minimal
$C(m)$ presentation, one may identify a redundant generator by seeking a
relation of the form $(a, w)$ or $(w,a)$ for some $a \in \scrA$ and $w \in \scrA^*$,
and then use the approach described above to effectively compute a
$C(m)$ presentation with strictly fewer generators. By iteration, one may
thus compute a generator-minimal $C(m)$ presentation for the same monoid.

Finally, to obtain a generator-minimal $C(m)$ equivalence presentation,
it suffices to take a generator-minimal $C(m)$ presentation and compute
the equivalence closure of the set of relations. By our remarks at the
beginning of the section, this is a $C(m)$ equivalence presentation for
the same monoid, and is clearly still generator-minimal since it has
the same generating set.
\end{proof}

\begin{corollary}\label{cor_mingens}
In any monoid admitting a $C(2)$ presentation, the set of non-identity
indecomposable
elements forms a generating set which is contained in every generating
set for the monoid.
\end{corollary}
\begin{proof}
By Proposition~\ref{prop_genminexists}, such a monoid has a generator-minimal
$C(2)$ presentation, that is, one in which no generators are redundant. By
Proposition~\ref{prop_indecomp}, the generators in this presentation must
be indecomposable, which shows that the set of all non-identity indecomposable
elements forms a generating set. Finally, we have already observed that
every non-identity indecomposable element
must lie in every generating set for the monoid.
\end{proof}

We now introduce some more terminology. Let $\presar$ and $\presbs$ be presentations. An
isomorphism $\hat\sigma : \scrA^* \to \scrB^*$ of free monoids is called an
\textit{inclusion} of $\presar$
into $\presbs$ if for every relation
$$(u, v) \in \scrR$$
there is a relation
$$(\hat\sigma(u), \hat\sigma(v)) \in \scrS.$$
We say that $\presar$ is a \textit{sub-presentation} of $\presbs$ if there
exists an inclusion of $\presar$ into $\presbs$. An \textit{isomorphism} between
two presentations is an inclusion whose inverse is also an inclusion, and
two presentations are called \textit{isomorphic} if there is an isomorphism
between them. Recall that, since a free monoid has a unique free generating
set, isomorphisms between free monoids are in a natural one-one correspondence
with bijections between the corresponding free generating sets. While
conceptually inclusions should be thought of as isomorphisms between free monoids,
computationally it is more practical to work with the corresponding bijections.
Notice that it is a completely routine matter to determine, given two
finite presentations and a bijection between the alphabets, whether the
bijection extends to an inclusion (and hence also whether it extends to
an isomorphism).

The following theorem says, informally, that any generator-minimal equivalence
presentation contains a copy of every generator-minimal $C(2)$
presentation for the same monoid.

\begin{theorem}\label{thm_inclusion}
Let $\presar$ be a generator-minimal $C(2)$ presentation and
$\presbs$ be any generator-minimal equivalence presentation. Let
$\sigma : \presar \to \presbs$ be an isomorphism between the monoids
presented. Then there is a unique morphism of free monoids
$\hat\sigma : \scrA^* \to \scrB^*$
such that $[\hat\sigma(a)]_\scrS = \sigma([a]_\scrR)$ for all $a \in \scrA$,
and this morphism is an inclusion of $\presar$ into $\presbs$.
\end{theorem}
\begin{proof}
Since $\presar$ is generator-minimal, the formal generators must
represent distinct elements of the monoid and, by Proposition~\ref{prop_indecomp},
those elements are indecomposable. Since $\sigma$ is an isomorphism, we
deduce that the images of these generators are indecomposable, which by
Corollary~\ref{cor_mingens} means that they lie in the set of elements
represented by the formal generators in $\scrB$. It follows that
we can define
$\hat\sigma : \scrA^* \to \scrB^*$ by chosing for each $a \in \scrA$ an
element $\hat\sigma(a) \in \scrB$ such that
$[\hat\sigma(a)]_\scrS = \sigma([a]_\scrR)$, and then extending to
a morphism of $\scrA^*$ using the universal property of free monoids.
Moreover, since $\presbs$ is
generator-minimal, the generators in $\presbs$ represent distinct
elements, so these choices are unique and $\hat\sigma$ is the unique
function with the given property.

Now $\sigma$ is injective, so $\hat\sigma$ separates $\scrA$, and
since $\scrB$ is free it follows that $\hat\sigma$ is injective.
Moreover, $\sigma$ is an isomorphism, so the images of the generating set for
$\presar$ comprise a generating set for $\presbs$. Since they are contained in
$\scrB$ and $\presbs$ is generator-minimal, this means that $\hat\sigma$
maps $\scrA$ bijectively onto $\scrB$, and hence is an isomorphism between
$\scrA^*$ and $\scrB^*$. This isomorphism clearly has the
property that $u \equiv_\scrR v$ if and only if $\hat\sigma(u) \equiv_\scrS \hat\sigma(v)$.

Now suppose $(x, y) \in \scrR$. To show that
$\hat\sigma$ is an inclusion, and hence complete the proof of the theorem, we
need to show that $(\hat\sigma(x), \hat\sigma(y)) \in \scrS$.
By the property of $\hat\sigma$ described above we certainly have that
$\hat\sigma(x) \equiv_\scrS \hat\sigma(y)$ in the monoid $\presbs$.
This means that we can choose a sequence of words
$$u_0 = \hat\sigma(x), \ u_1, \ u_2, \ \dots, \ u_n = \hat\sigma(y) \ \in \ \scrB^*$$
such that each $u_{i+1}$ can be obtained from $u_i$ by an application of a
relation from $\scrS$. This by definition means that for $1 \leq i < n$ we may choose
$p_i, q_i, q_i', r_i \in \scrB^*$ such that $u_i = p_i q_i r_i$,
$u_{i+1} = p_i q_i' r_i$, $q_i \neq q_i'$ and (using the fact that
$\scrS$ is symmetric) also $(q_i, q_i') \in \scrS$.

For each $i$, let $w_i = \hat\sigma^{-1}(u_i) \in \scrA^*$. Then
using again the property above of $\hat\sigma$ we have that
$w_1, \dots, w_n$
are all equivalent in $\presar$ to the relation word $w_0 = \hat\sigma^{-1}(u_0) = x$.
By Proposition~\ref{prop_relationwordclass}, this
means in particular that every $w_i$ is a relation word in the presentation
$\presar$.

We claim now that $u_i = q_i$ for every $i$. Indeed, suppose for a contradiction
that there is some $i$ for which $q_i \neq u_i$. 
We know that $q_i \equiv_\scrS q_i'$ and
$q_i \neq q_i'$ (since the rewrite is non-trivial), which by the property of
$\hat\sigma$ noted above means that
$\hat\sigma^{-1}(q_i) \equiv_\scrR \hat\sigma^{-1}(q_i')$
and $\hat\sigma^{-1}(q_i) \neq \hat\sigma^{-1}(q_i')$.
But $\hat\sigma^{-1}(q_i)$ is a proper
factor of $w_i = \hat\sigma^{-1}(u_i)$, which is a relation word, so
by Proposition~\ref{prop_relationwordfactorclass} we must have
$\hat\sigma^{-1}(q_i) = \hat\sigma^{-1}(q_i')$ and hence $q_i = q_i'$
which gives the desired contradiction and completes the proof of the
claim that $u_i = q_i$ for every $i$.

It follows immediately that $q_i' = u_{i+1} = q_{i+1}$.
But now $(u_i, u_{i+1}) = (q_i, q_i') \in \scrS$ for every $i$,
which since $\scrS$ is transitive, means that
$(u_0, u_n) \in \scrS$ as required.
\end{proof}

\begin{corollary}\label{cor_isocondition}
Two generator-minimal $C(2)$ equivalence presentations present isomorphic
monoids if and only if they are isomorphic.
\end{corollary}
\begin{proof}
It is obvious that isomorphic presentations present isomorphic monoids.

Conversely, let $\presar$ and $\presbs$ be generator-minimal
$C(2)$ equivalence presentations, and let $\sigma_1 : \presar \to \presbs$ be an
isomorphism between the monoids presented, and $\sigma_2 : \presbs \to \presar$
be its inverse.
Let $\hat\sigma_1 : \scrA^* \to \scrB^*$ and $\hat\sigma_2 : \scrB^* \to \scrA^*$
be the inclusions given by Theorem~\ref{thm_inclusion}. To show that the
presentations are isomorphic, it will suffice to show that the inclusions
$\hat\sigma_1$ and $\hat\sigma_2$ are mutually inverse. Now we know from
Theorem~\ref{thm_inclusion} that
$\hat\sigma_2$ is the unique isomorphism from $\scrB^*$ to $\scrA^*$ satisfying
$$[\hat\sigma_2(b)]_\scrR = \sigma_2([b]_\scrS)$$
for all $b \in \scrB$ so it will suffice to show that $\hat\sigma_1^{-1}$ also
satisfies this condition.

By the properties of $\hat\sigma_1$ guaranteed by
Theorem~\ref{thm_inclusion} we have
$$[\hat\sigma_1(a)]_\scrS = \sigma_1([a]_\scrR)$$
for all $a \in \scrA$. Applying $\sigma_1^{-1} = \sigma_2$ to both sides
we obtain
$$\sigma_2 ([\hat\sigma_1(a)]_\scrS) = \sigma_2 (\sigma_1 ([a]_\scrR)) = [a]_\scrR.$$
Now for any $b \in \scrB$, letting $a = \hat\sigma_1^{-1}(b)$ in this
expression yields
$$\sigma_2 ([b]_\scrS) = [\hat\sigma_1^{-1}(b)]_\scrR$$
as required.
\end{proof}

As a straightforward consequence, we obtain the fact that the isomorphism
problem for $C(2)$ monoids is algorithmically solvable.

\begin{theorem}\label{thm_isomorphism}
There is an algorithm which, given as input two $C(2)$ presentations,
decides whether the monoids presented are isomorphic.
\end{theorem}
\begin{proof}
Given two $C(2)$ presentations, by Proposition~\ref{prop_genminexists},
we may compute generator-minimal $C(2)$ equivalence presentations for the
same monoids. Now by Corollary~\ref{cor_isocondition} it suffices to check
if the resulting presentations are themselves isomorphic. This can clearly
be done, for example by enumerating all bijections between the generating
sets and checking if any extend to isomorphisms.
\end{proof}

The following theorem says that a small overlap presentation
without redundant generators will always satisfy the strongest small overlap
conditions satisfied by any presentation for the same monoid.

\begin{theorem}\label{thm_twoimpliesall}
Let $m$ be a positive integer, and $\presar$ be a generator-minimal
$C(2)$ presentation for a monoid which admits a $C(m)$ presentation.
Then $\presar$ is a $C(m)$ presentation.
\end{theorem}
\begin{proof}
If $m \leq 2$ then the claim is trivial, so suppose $m \geq 3$.

Suppose for a contradiction that $\presar$ does not satisfy $C(m)$. Then
there is a relation word $r \in \scrA^*$ which can be written as a product
of strictly fewer than $m$ pieces of $\presar$, say $r = r_1 r_2 \dots r_j$
where $r_1, \dots, r_j \in \scrA^*$ are pieces and $j < m$.

Now the monoid presented by $\presar$ admits a $C(m)$ presentation, so
by Proposition~\ref{prop_genminexists} it admits a generator-minimal
$C(m)$ equivalence presentation. Let $\presbs$ be such a presentation
and $\sigma : \presar \to \presbs$ an isomorphism. Let $\hat\sigma : \scrA^* \to \scrB^*$ be
the corresponding morphism of free monoids given by Theorem~\ref{thm_inclusion}.
Since $\hat\sigma$
is an inclusion it follows easily that $\hat\sigma(r)$ is
a relation word in $\presbs$, and $\hat\sigma(r_1), \dots, \hat\sigma(r_j)$
are pieces of $\presbs$. But since $\hat\sigma$ is a morphism, we have
$$\hat\sigma(r) = \hat\sigma(r_1) \hat\sigma(r_2) \dots \hat\sigma(r_j)$$
which contradicts the assumption that $\presbs$ satisfies the condition
$C(m)$.
\end{proof}

We also obtain a corresponding statement for strong small overlap
conditions.

\begin{corollary}\label{cor_twoimpliesallstrong}
Let $m$ be a positive integer, and $\presar$ be a generator-minimal
strongly $C(2)$ presentation for a monoid which admits a strongly
$C(m)$ presentation. Then $\presar$ is a strongly $C(m)$ presentation.
\end{corollary}
\begin{proof}
By Theorem~\ref{cor_twoimpliesallstrong} $\presar$ is a $C(m)$
presentation. Since it is strongly $C(2)$ it contains no repeated
relation words, which means it is strongly $C(m)$.
\end{proof}

\begin{corollary}\label{cor_compute}
There is an algorithm which, given as input a monoid presentation,
finds a generator-minimal $C(m)$ equivalence presentation for the same
monoid with $m \geq 2$ as high as possible, provided such a presentation
exists (and may not terminate otherwise).
\end{corollary}
\begin{proof}
If the presentation satisfies $C(m)$ for some $m \geq 2$ then by definition
it satisfies $C(2)$, and by Proposition~\ref{prop_genminexists} the
monoid presented admits
a generator-minimal $C(2)$ equivalence presentation. By
Theorem~\ref{thm_twoimpliesall}, any such presentation must also satisfy
$C(m)$, so it suffices to find such a presentation. It is well known and
easy to show (using for example the theory of Tietze transformations
\cite[Section~7.2]{Book93}) that there is an algorithm which,
given as input any monoid presentation, recursively enumerates all
presentations isomorphic to it. So it suffices to do this, checking
if each is a generator-minimal $C(2)$ equivalence presentation, until
we find one which is.
\end{proof}

\begin{corollary}
Let $m$ be a positive integer. Then there is no algorithm to decide,
given as input a monoid presentation, whether the monoid presented admits
a $C(m)$ presentation.
\end{corollary}
\begin{proof}
It is well known (see for example \cite[Corollary~7.3.8]{Book93}) that there
is no algorithm to decide, given a monoid presentation, whether the monoid presented is trivial. 
The trivial monoid admits the empty presentation, which is a $C(m)$
presentation for every $m$. If there were an algorithm to decide whether
a given monoid admits a $C(m)$ presentation then one could decide whether
the monoid was trivial, by first checking if it admits a $C(m)$ presentation.
If not then it cannot be trivial. If so then by Corollary~\ref{cor_compute}
we could compute a $C(m)$ presentation for the monoid, and by
Theorem~\ref{thm_isomorphism} we could check if it is isomorphic to the
trivial monoid.
\end{proof}

Using Theorem~\ref{thm_isomorphism} and techniques from \cite{K_generic},
we can obtain asymptotics for the number of non-isomorphic semigroups
admitting $a$-generator, $k$-relation presentations of a certain length.
Recall that for functions $f, g : \mathbb{N} \to \mathbb{N}$ we say that
$f(n)$ is $O(g(n))$ [respectively, $\Omega(g(n)$, $\theta(g(n))$] if $f(n)$
is bounded above [respectively, below, above and below] by a linear function
of $g(n)$ with positive coefficients.

\begin{theorem}\label{thm_noniso}
Let $a$ and $k$ be fixed positive integers. Then there are $\theta(a^n n^{2k-1})$
distinct isomorphism types of semigroups admitting $a$-generator, $k$-relation
presentations of sum relation length $n$.
\end{theorem}
\begin{proof}
The proof utilises some arguments from \cite{K_generic}; since this
result is not the main purpose of the present paper we refrain from repeating
these in detail and instead refer the reader to \cite{K_generic} for a
more detailed explanation whether appropriate.

For ease of counting, we will first consider \textit{ordered
presentations}, by which we mean a pairs $\langle \scrA \mid \scrR \rangle$
where $\scrA$ is an alphabet and $\scrR$ is a (finite or infinite)
\textit{sequence} of pairs of words over $\scrA$. As observed in \cite{K_generic},
an ordered $\scrA$-generated, $k$-relation monoid presentation of sum relation length
$n$ is uniquely determined by its \textit{shape} (the sequence of lengths
of relation words) together with the concatenation in order of the relation
words (a word of length $n$).

Recall (from for example
 \cite[Theorem~5.2]{Bona02}) that the number of weak compositions of
$s$ into $r$ (that is, ordered sequences of $r$
non-negative integers summing to $s$) is given by
$$C'_r(s) = \frac{(s+r-1)!}{s! (r-1)!}.$$
Note that if $r$ is fixed and $s$ remains variable then $C'_r(s)$ is a
polynomial of degree $r-1$ in $s$.

Let $f: \mathbb{N} \to \mathbb{R}$ be a function.
The total number of shapes of length $n$ is $C'_{2k}(n)$. As shown in the
proof of \cite[Lemma~3.4]{K_generic}, the number of shapes featuring a
block of size $f(n)$ or less is bounded above by
$$2k (f(n)+1) C'_{2k-1}(n)$$
so the number of such shapes \textit{not} featuring such a block is
bounded below by 
\begin{align*}
C'_{2k}(n) - 2k (f(n)+1) C'_{2k-1}(n) &= \frac{(n+2k-1)!}{n!(2k-1)!} - 2k (f(n)+1) \frac{(n+2k-2)!}{n! (2k-2)!} \\
&= \frac{(n+2k-2)!}{n!(2k-2)!} \left( \frac{n+2k-1}{2k-1} - 2k (f(n)+1) \right) \\
&= C'_{2k-1}(n) \left( \frac{n}{2k-1} - 2k (f(n) + 1) + 1 \right) \\
&= n C'_{2k-1}(n) \left( \frac{1}{2k-1} - \frac{2k (f(n) + 1) + 1}{n} \right)
\end{align*}

Now fix an alphabet $\scrA$ of size $a$.
By \cite[Lemma~3.2]{K_generic} there are at most $n^2 a^{n-f(n)}$ distinct
words of length $n$ which contain a repeated factor of length $f(n)$, and
so there are at least
$$a^n - n^2 a^{n-f(n)} = a^n \left(1-\frac{n^2}{a^{f(n)}} \right)$$
words of length $n$ which do not contain such a repeated factor. Combining
such a word with a shape with no blocks of length $f(n)$ or less clearly
yields a presentation in which no relation word appears as a factor of any
other relation word, that is a strongly $C(2)$ presentation. Dividing by
$k!$ to allow for reordering the $k$ relations,
we have at least
$$\frac{1}{k!} \ a^n \ \left(1-\frac{n^2}{a^{f(n)}} \right) \ C'_{2k-1}(n) \ n \ \left( \frac{1}{2k-1} - \frac{2k (f(n) + 1) + 1}{n} \right)$$
$\scrA$-generated $k$-relation generator-minimal strongly $C(2)$ unordered
presentations of length $n$. 

Since a presentation presents the same monoid as its equivalence closure
it follows by Corollary~\ref{cor_isocondition} that two such presentations
present isomorphic monoids if and only their equivalence closures are
isomorphic presentations. Since a strongly $C(2)$ presentation must be transitive and
irreflexive, the only difference between two strongly $C(2)$ presentations with
the same equivalence closure is in the order of each relation, so at
most $2^k$ distinct strongly $C(2)$ presentations have the same equivalence
closure. And since an isomorphism of presentations over $\scrA$ is determined by
a bijection on $\scrA$, each such equivalence closure is isomorphic
to at most $a!$ others. Thus, at most $2^k a!$ distinct strongly $C(2)$
presentations present isomorphic monoids, so the strongly $C(2)$ presentations found
above must present at least
$$\frac{1}{a! k! 2^k} \ a^n \ \left(1-\frac{n^2}{a^{f(n)}} \right) \ C'_{2k-1}(n) \ n \ \left( \frac{1}{2k-1} - \frac{2k (f(n) + 1) + 1}{n} \right)$$
distinct isomorphism types of monoids. Setting $f(n) = 3 \log_a n$ this becomes
at least
$$\frac{1}{a! k! 2^k} \ a^n \ \left(1-\frac{1}{n} \right) \ C'_{2k-1}(n) \ n \ \left( \frac{1}{2k-1} - \frac{2k (3 \log_a n + 1) + 1}{n} \right)$$
distinct isomorphism types.
Now the factor
$$\frac{1}{a! k! 2^k}$$
is a positive constant, while the factors
$$\left(1-\frac{1}{n} \right) \text{ and } \left( \frac{1}{2k-1} - \frac{2k (3 \log_a n + 1) + 1}{n} \right)$$
are eventually bounded below by positive constants
and $C'_{2k-1}(n)$ is a polynomial of degree $2k-2$ in $n$. Thus, the number
of distinct isomorphism types is $\Omega(a^n n^{2k-1})$.

Finally, the number of isomorphism types of semigroups admitting $a$-generator,
$k$-relation presentations of length $n$ is clearly bounded above by the number of
$k$-generator presentations over a fixed alphabet of size $a$. This is bounded
above by the number of words of length $n$ times the number of weak compositions
of $n$ into $2k$, that is $a^n C'_{2k}(n)$. Since $C'_{2k}(n)$ is a 
polynomial of degree $2k-1$ in $n$, it follows that the number of
presentations, and hence the number of isomorphism types, is $O(a^n n^{2k-1})$.
\end{proof}

A closer analysis of the combinatorics in the proof of Theorem~\ref{thm_noniso}
would of course yield more precise information about the number of isomorphism
classes, as well as asymptotics applicable when $a$ and $k$ are permitted to
vary.

\section{Cancellativity}\label{sec_canc}

In this section we investigate the conditions under which the monoid
presented by a $C(4)$ presentation is left cancellative, right
cancellative or (two-sided) cancellative. It transpires that these
properties can be characterised by very simple and natural conditions
on the presentation, and from this it follows that one can check in
linear time whether a given $C(4)$ monoid has any of these properties.
Interestingly, it follows also that cancellativity properties are
\textit{asymptotically visible} properties of finite monoid presentations,
in the sense that the probability that an $\scrA$-generated, $k$-relation
presentation of length selected uniformly at random presents a (left,
right or two-sided) cancellative monoid
converges to a value strictly between $0$ and $1$, as the size of the
presentation increases.

We begin with a technical definition. Let $\presar$ be a $C(4)$
presentation. We define a function $\rho : \scrA^* \to \mathbb{Z}$ as
follows. For $w \in \scrA^*$ define $\rho(w) = -1$ if $w$ has no clean
overlap prefix; otherwise we define $\rho(w)$ to be the length of the
suffix of $w$ following the clean overlap prefix. Notice that $\rho$
takes integer values greater than or equal to $-1$; we shall use it
as an induction parameter.

\begin{lemma}\label{lemma_rdrops}
If $w = XYZw'$ where $XY$ is a clean overlap prefix and $p$ is a piece
then $\rho(p w') < \rho(w)$ 
\end{lemma}
\begin{proof}
By definition we have $\rho(w) = |Zw'| \geq 0$. If $p w'$ has no clean
overlap prefix then again by definition we have $\rho(pw') = -1$ and we are
done. Otherwise, suppose $pw'$ has a clean overlap prefix $aX'Y'$, say
$pw' = aX'Y'w''$. Then $\rho(pw') = |w''|$.
Now $X'$ is the maximum piece prefix of a relation word $X'Y'Z'$, so
$X'Y'$ is not a piece. In particular $X'Y'$ cannot be a factor of $p$,
so we must have that $aX'Y'$ is strictly longer than $p$, and so $|w''| < |w'|$.
Hence,
$$\rho(pw') = |w''| < |w'| \leq |Zw'| = \rho(w)$$
as required.
\end{proof}

We are now ready to prove our main result of this section, which
characterises $C(4)$ equivalence presentations of left cancellative monoids.

\begin{theorem}\label{thm_canc}
Let $\presar$ be a $C(4)$ equivalence presentation. Then the monoid
presented is left cancellative if and only if $\scrR$ contains no relation of
the form $(ar, as)$ where $a \in \scrA$ and $r,s \in \scrA^*$ with $r \neq s$. 
\end{theorem}
\begin{proof}
Suppose first that $\scrR$ contains a relation $(ar, as)$ for some
letter $a \in \scrA$ and words $r, s \in \scrA^*$ with $r \neq s$. Since
the presentation satisfies $C(2)$ and $r$ is a proper factor of a relation
word, it follows by Proposition~\ref{prop_relationwordfactorclass} that $r \nequiv s$.
But by definition we have $ar \equiv as$, so the monoid presented is not
left cancellative.

Conversely, suppose $\presar$ contains no relations of the given form, and
suppose for a contradiction that the monoid presented by $\presar$ is not
left cancellative. Then the generators are not left
cancellable, so there exists a generator $a \in \scrA$ and words $u,v \in \scrA^*$
such that $au \equiv av$ but $u \nequiv v$. Let $a$, $u$ and $v$ be chosen
such that these conditions hold and $\rho(au)$ is as small as possible.

First, note that if $au$ has no clean overlap prefix then by
Proposition~\ref{prop_nomopnorel} we have $au = av$ as words, so that
$u = v$ which gives a contradiction. Next, suppose that $au$ has a
clean overlap prefix $wXY$ with $w$ non-empty, say $au = wXYu''$. Then
by Proposition~\ref{prop_dumpprefix} we have $av = wv'$ where
$v' \equiv XYu''$. But now $w = aw'$ for some word $w$ and we have
$u = w'XYu'' \equiv w'v' = v$, again giving a contradiction.

There remains only the case that $au$ has a clean overlap prefix of the
form $XY$. Then one of the six mutually exclusive conditions of
Lemma~\ref{lemma_eq} hold in respect of $au$ and $av$.
Now conditions (3)-(6) are impossible because $au$ and $av$ begin with
the same letter, while $XY$ and the complement $\ol{XY}$ cannot.
Indeed, if they did begin with the same letter then so would the distinct
relation words $XYZ$ and $\ol{XYZ}$. But by symmetry and transitivity of
the presentation we have $(XYZ, \ol{XYZ}) \in \scrR$, so this would contradict the assumption on
the presentation.

If condition (2) holds we have $au = XYu'$, $av = XYv'$ where $u' \equiv v'$.
But $Y$ is non-empty, so we may write $XY = ar$ for some $r \in \scrA^*$,
whereupon $u = ru' \equiv rv' = v$ again giving a contradiction.

If condition (1) holds then $au = XYZu''$ and $av=XYZv''$ and
$\ol{Z} u'' \equiv \ol{Z} v''$ for some complement $\ol{Z}$ of $Z$.
We claim that the existence of any piece $q$ such that $qu'' \equiv qv''$
implies that $u'' \equiv v''$; the proof of this claim is by induction on
the length of $q$. For the base case, if $q$ is the empty word then certainly we have
$u'' = qu'' \equiv qv'' = v''$. Now suppose for induction that $q$ is
not the empty word, that $qu'' \equiv qv''$ and that the claim holds for
shorter $q$. Write $q = bq'$ where $b \in \scrA$ and $q' \in \scrA^*$.
Then by Lemma~\ref{lemma_rdrops} we have
$$\rho(bq'u'') \ = \ \rho(qu'') \ < \ \rho(XYZu'') \ = \ \rho(au)$$
so by the minimality assumption on $a$, $u$ and $v$ we have
$q'u'' \equiv q'v''$, which by the inductive hypothesis implies that
$u'' \equiv v''$. This completes the proof of the claim.
Now returning to the main proof, setting $q = \ol{Z}$ we deduce that
$u'' \equiv v''$, and an argument similar to that in case (2) completes
the proof that $u \equiv v$, once again establishing the required
contradiction.
\end{proof}

The combinatorial condition given by Theorem~\ref{thm_canc} can
clearly be checked in polynomial time (more precisely, linear time
in the RAM model of computation) so we have:

\begin{corollary}
There is an algorithm which, given as input a $C(4)$ equivalence presentation,
decides in polynomial time whether the monoid presented is left cancellative,
right cancellative and/or cancellative.
\end{corollary}

\begin{corollary}\label{cor_strongcanc}
Let $\presar$ be a strongly $C(4)$ presentation. Then the monoid
presented is left cancellative if and only if $\scrR$ contains no relation of
the form $(ar, as)$ where $a \in \scrA$ and $r,s \in \scrA^*$. 
\end{corollary}
\begin{proof}
Since a strongly $C(4)$ presentation contains no repeated relation words
it is already transitive. Hence its equivalence closure $\presas$ can be
obtained simply by adding all relations of the form $(u,u)$, $(v,v)$ and
$(v,u)$ where $(u,v)$ is a relation in $\scrR$. Now $\presas$ presents the
same monoid as $\presar$, so by Theorem~\ref{thm_canc}, $\presar$ is
left cancellative exactly if $\scrS$ contains no relation of the form
$(ar, as)$ with $a \in \scrA$ and $r \neq s$. Clearly, this is true
exactly if $\scrR$ contains no relation of the form $(ar,as)$ with
$a \in \scrA$ and $r,s \in \scrA^*$.
\end{proof}

Our results in this section also have a consequence for the theory of
generic properties in monoids and semigroups \cite{K_generic}.
In the
terminology of generic-case complexity, the following theorem says that
left, right and two-sided cancellativity
are \textit{aymptotically visible} properties of $A$-generated, $k$-relation
monoids. Such properties are of interest since they seem to be rather rare,
with most abstract properties of groups and semigroups tending to be either
\textit{generic} or \textit{negligible} amongst finitely presented examples.
To avoid defining large amounts of terminology for a single use, we
state the result in elementary combinatorial terms.

\begin{theorem}\label{thm_cancprop}
Let $\scrA$ be an alphabet with $|\scrA| \geq 2$. Then the proportion $\scrA$-generated,
$k$-relation monoid presentations of sum relation length $n$ (or of maximum
relation length $n$) which present left cancellative [right
cancellative] monoids approaches
$$\left( \frac{|\scrA|-1}{|\scrA|} \right)^{k}.$$
as $n$ tends in $\infty$. The proportion of $\scrA$-generated,
$k$-relation monoid presentations of sum relation length $n$ (or of maximum
relation length $n$) which present
cancellative monoids approaches
$$\left( \frac{|\scrA|-1}{|\scrA|} \right)^{2k}.$$
as $n$ tends in $\infty$. In particular, for $|\scrA| \geq 2$, left [right]
cancellativity is an asymptotically visible property of $A$-generated
$k$-relation presentations.
\end{theorem}
\begin{proof}
As with the proof of Theorem~\ref{thm_noniso}, we make use of some
results and arguments from \cite{K_generic}. Rather than repeating these at
length, we instead refer the interested reader to that paper.

By \cite[Theorem~3.5]{K_generic} (or \cite[Theorem~3.8]{K_generic} when the maximum
relation length is considered), the proportion of ordered presentations which
are not strongly $C(4)$ approaches $0$ as $n$ increases. Corollary~\ref{cor_strongcanc} tells
us that a strongly $C(4)$ presentation presents a left [right] cancellative monoid
exactly if it contains no relation of the form $ar = as$ [$ra = sa$] where $a \in \scrA$ and
$r, s \in \scrA^*$, so to find the limit of the proportion of left [right]
cancellative monoid, it suffices to compute the proportion of presentations
satisfying this condition.

As described in the proof of Theorem~\ref{thm_noniso},
an ordered presentation is uniquely
determined by its shape (a weak composition of $n$ into $2k$) and the
concatenation in order of its relation words (a word over $\scrA$ of length
$n$). Now for any shape not featuring a block of length $0$, the proportion
of words which do not yield two relation words beginning [ending] with the
same letter is clearly exactly $( (|\scrA|-1) / |\scrA|)^k$. 
Similarly, for any shape not featuring a block of length $0$ or $1$ (so
that the first and last positions of each relation word are distinct), the
proportion of words which do not yield two relation words beginning with the same
letter or ending with the same letter is exactly
$( (|\scrA|-1) / |\scrA|)^{2k}$.

Now by \cite[Lemma~3.4]{K_generic} (or \cite[Lemma~3.6]{K_generic} when the
maximum  relation length is considered), the proportion of $\scrA$-generated,
$k$-relation ordered presentations of length $n$ whose shape features a relation word
of length $0$ or $1$ tends to zero as $n$ tends to infinity. It follows that
we may ignore presentations with these shapes, and conclude that the
proportion of $\scrA$-generated, $k$-relation ordered presentations of sum or relation
length $n$ which are left cancellative, right cancellative and cancellative
approaches
$$\left( \frac{|\scrA|-1}{|\scrA|} \right)^k, \ \left( \frac{|\scrA|-1}{|\scrA|} \right)^k \text{ and } \left( \frac{|\scrA|-1}{|\scrA|} \right)^{2k}$$
respectively as $n$ tends to $\infty$.

Finally, we turn our attention to unordered presentations. By
\cite[Theorem~4.10]{K_generic} the proportion of unordered presentations
of length $n$ which fail to be strongly $C(4)$ approaches $0$ at $n$ tends to 
$\infty$, so it suffices to compute the proportion of strongly $C(4)$
unordered presentations which present left cancellative [right cancellative,
cancellative] monoids. However,
since a strongly $C(4)$ presentation contains no repeated relation words, each
$k$-relation unordered strongly $C(4)$ presentation corresponds to precisely
$k!$ distinct ordered strongly $C(4)$ presentations. Since moreover
every ordered strongly $C(4)$ presentation arises in this way, it follows
that the required proportions are the same for unordered presentations as for
ordered presentations.
\end{proof}

In fact the argument in the last paragraph of the proof of
Theorem~\ref{thm_cancprop} establishes the following fact, which
generalises \cite[Theorem~4.10]{K_generic}.
\begin{theorem}
Let $\mathscr{C}$ be any class of abstract monoids, $\scrA$ an alphabet with
$|\scrA| \geq 2$ and $k$ a positive integer. Then the proportion of
$\scrA$-generated $k$-relation ordered monoid presentations of sum
[maximum] relation length $n$ which present monoids in $\mathscr{C}$
converges with the proportion of
$\scrA$-generated $k$-relation monoid presentations of sum
[maximum] relation length $n$
which present monoids in $\mathscr{C}$ as $n$ tends to $\infty$.
\end{theorem}

\section*{Acknowledgements}

This research was supported by an RCUK Academic Fellowship. The author
would like to thank S.~W.~Margolis and N.~Ruskuc for posing some of
the questions considered in this paper.

\bibliographystyle{plain}

\def\cprime{$'$} \def\cprime{$'$}

\end{document}